\documentclass[a4paper]{amsart}
\usepackage{amsmath} 
\usepackage{amssymb}
\usepackage{latexsym}
\usepackage[mathscr]{eucal}

\font\est=eufb10 scaled\magstep1
\font\es=eufm10

\newcommand{\so}{\mathfrak{so}}
\newcommand{\gge}{\mathfrak{e}}
\newcommand{\gP}{\mathfrak{P}}
\newcommand{\gW}{\mathfrak{W}}
\newcommand{\gJ}{\mathfrak{J}}
\newcommand{\gC}{\mathfrak{C}}
\newcommand{\spin}{\mathfrak{spin}}

\def\gge{\mbox{\es{e}}} 
\def\gf{\mbox{\es {f}}}

\def\gC{\mbox{\es {C}}}

\def\gJ{\mbox{\es {J}}}

\def\gP{\mbox{\es {P}}}

\def\gW{\mbox{\es {W}}}

\def\sl{\mbox{\es {sl}}}

\def\so{\mbox{\es {so}}}

\def\spin{\mbox{\es {spin}}}
\def\dfrac#1#2{\displaystyle \frac{#1}{#2}}

\def\det{\mbox{\rm {det}}}
\def\diag{\mbox{\rm {diag}}}

\def\dlim{\displaystyle \lim}

\def\ad{\mbox{\rm {ad}}}

\def\Hom{\mbox{\rm {Hom}}}

\def\Iso{\mbox{\rm {Iso}}}

\def\Ker{\mbox{\rm {Ker}}}

\def\ov{\overline}
\def\wti{\widetilde}

\def\dsum{\displaystyle \sum}

\def\dfrac#1#2{\displaystyle \frac{#1}{#2}}

\def\R{\mbox{\boldmath $R$}}

\def\Z{\mbox{\boldmath $Z$}}

\def\0{\mbox{\boldmath {0}}}    
\def\1{\mbox{\boldmath {1}}}      
\def\2{\mbox{\boldmath {2}}}      
\def\3{\mbox{\boldmath {3}}}      
\def\4{\mbox{\boldmath {4}}}      
\def\5{\mbox{\boldmath {5}}}      
\def\6{\mbox{\boldmath {6}}}      
\def\7{\mbox{\boldmath {7}}}      
\def\8{\mbox{\boldmath {8}}}      
\def\9{\mbox{\boldmath {9}}}

\def\gge{\mbox{\est {e}}}

\def\sl{\mbox{\es {sl}}} 
\def\so{\mbox{\es {so}}}

\usepackage{epic,eepic}

\title[Fixed points subgroups $G^{\sigma,\sigma'}$ by two involutive automorphisms 
$\sigma$, $\sigma'$]
{Fixed points subgroups $G^{\sigma,\sigma'}$ by two involutive automorphisms 
$\sigma$, $\sigma'$ of exceptional compact Lie group $G$, Part II, $G = E_8$}

\dedicatory{Dedicated to Professor Ichiro Yokota on the occasion of his eightieth birthday}
 
\author[Toshikazu Miyashita]{Toshikazu Miyashita}

\begin{document}

\makeatletter
\begin{abstract}
For the simply connected compact exceptional Lie group $E_8$, we determine the structure of subgroup $(E_8)^{\sigma, \sigma'}$ of $E_8$ which is the intersection $(E_8)^\sigma \cap (E_8)^{\sigma'}$. \vspace{0.5mm}Then the space $E_8/(E_8)^{\sigma, \sigma'}$ is the exceptional $\mathbb{Z}_2 \times \mathbb{Z}_2$- symmetric space of type 
$\text{E\@Roman{8}-\@Roman{8}-\@Roman{8}}$, and that we give two involutions $\sigma, \sigma'$ for the space $E_8/(E_8)^{\sigma, \sigma'}$ concretely.
\end{abstract}
\maketitle

\vspace{-8mm}

\section{\bf Introduction}
\vspace{3mm}

In the preceding  paper [3], for the simply connected compact exceptional Lie 
groups $G = F_4, E_6$ and $E_7$, we determined the structure of subgroups 
$G^{\sigma, \sigma'}$ of $G$ which is the intersection $G^\sigma \cap G^{\sigma'}$. Their results were given as
\begin{eqnarray*}
        ({F_4})^{\sigma,\sigma'} \!\!\! & \cong &\!\!\! Spin(8),
\\[0.5mm]
        ({E_6})^{\sigma,\sigma'} \!\!\! &\cong &\!\!\! (U(1) \times U(1) \times Spin(8))/(\Z_2 
\times \Z_4),
\vspace{1mm}\\
        ({E_7})^{\sigma,\sigma'} \!\!\! &\cong &\!\!\! (SU(2) \times SU(2) \times SU(2) \times 
Spin(8))/(\Z_2 \times \Z_2).
\end{eqnarray*}

\noindent In [1], a classification was given of the exceptional $\mathbb{Z}_2 \times \mathbb{Z}_2$-symmetric spaces $G/K$ where $G$ is an exceptional compact Lie groups and $Spin(8)$, and the structure of $K$ was determined as Lie algebras. 
Our results correspond to the realizations of the exceptional $\mathbb{Z}_2 \times \mathbb{Z}_2$-symmetric spaces with the pair $(\sigma,\sigma')$ of commuting involutions of $F_4, E_6$ and $E_7$.
To be specific, since the groups $(F_4)^{\sigma, \sigma'},\, (E_6)^{\sigma, \sigma'}$ and $\,(E_7)^{\sigma, \sigma'}$ are connected from the results above, we see that the space $F_4/(F_4)^{\sigma, \sigma'}, \,$ $E_6/(E_6)^{\sigma, \sigma'}$ and $E_7/(E_7)^{\sigma, \sigma'}$ are the exceptional $\mathbb{Z}_2 \times \mathbb{Z}_2$- symmetric spaces of type 
$\text{F\@Roman{2}-\@Roman{2}-\@Roman{2}, E\@Roman{3}-\@Roman{3}-\@Roman{3}\,\,and\,\,E\@Roman{6}-\@Roman{6}-\@Roman{6}}$, respectively.

In the present paper, for the simply connected compact exceptional Lie group $E_8$, 
we consider the subgroup $(E_8)^{\sigma,\sigma'}$ of ${E_8}$. We showed in [4] that the Lie algebra $(\mathfrak{e}_8)^{\sigma,\sigma'}$ of $(E_8)^{\sigma,\sigma'}$ is isomorphic to $\mathfrak{so}(8) \oplus \mathfrak{so}(8)$ as a Lie algebra. So, it is main purpose to determine the structure of the subgroup $({E_8})^{\sigma,\sigma'}$\!\!\!. Our result is as follows:
$$
     (E_8)^{\sigma, \sigma'} \cong (Spin(8) \times Spin(8))/(\Z_2 \times \Z_2). 
$$
The first $Spin(8)$ is $({E_8})^{\sigma,\sigma',\mathfrak{so}(8)} \subset (E_8)^{\sigma, \sigma'}$ and the second $Spin(8)$ is $(F_4)^{\sigma,\sigma'} \subset (E_8)^{\sigma, \sigma'}$, where the definition of the group $({E_8})^{\sigma,\sigma',\mathfrak{so}(8)}$ is showed in section 4.
On the exceptional $\mathbb{Z}_2 \times \mathbb{Z}_2$-symmetric spaces for 
$E_8$, there exist four types ([1]). This amounts to the realization of one type of these four types with the pair $(\sigma, \sigma')$ of commuting involutions of $E_8$.    

The essential part to prove this fact is to show the connectedness of the group $(E_8)^{\sigma,\sigma',\mathfrak{so}(8)}$. For this end, we need to treat the complex case, that is, we need the following facts.
\begin{eqnarray*}
      ({F_4}^C)^{\sigma,\sigma'}\!\!\! &\cong & \!\!\! Spin(8, C),
\\[0.5mm]
      ({E_7}^C)^{\sigma,\sigma',\mathfrak{so}(8)}\!\!\! &\cong & \!\!\! SL(2, C) \times SL(2, C) \times SL(2, C), 
\\[0.5mm]
      ({E_7}^C)^{\sigma,\sigma'}\!\!\! &\cong & \!\!\! (SL(2, C) \times SL(2, C) \times SL(2, C) \times Spin(8, C))/(\Z_2 \times \Z_2),
\\[0.5mm]
 \hspace*{-10mm} ({E_8}^C)^{\sigma,\sigma',\mathfrak{so}(8)}\!\!\! &\cong & \!\!\! Spin(8. C),
\\[0.5mm]
      ({E_8}^C)^{\sigma,\sigma'}\!\!\! &\cong & \!\!\!(Spin(8, C) \times Spin(8, C))/(\Z_2 \times \Z_2)
\end{eqnarray*}

and the connectedness of $({E_8}^C)^{\sigma,\sigma', \mathfrak{so}(8, C)}$. Even if some of their proofs are similar to the previous paper [5] and [6], we write in detail again.
\vspace{1mm}

To consider the group $E_8$, we need to know some knowledge of the group $E_7$. As for them we refer [6] and [7]. This paper is a continuation of [3]. We use the same notations as [3],[6] and [7]. 
\vspace{2mm}

\section
{\bf Lie groups ${F_4}^C, F_4, {E_6}^C, {E_7}^C, {E_8}^C$ and $E_8$}
\vspace{2mm}

We describe definitions of Lie groups used in this paper.
\vspace{1mm}

Let $\gJ^C$ and $\gJ$ be the exceptional $C$- and $\R$-Jordan algebras, respectively. The connected complex Lie group ${F_4}^C$ and the connected compact Lie group $F_4$ are defined by
\begin{eqnarray*}
      {F_4}^C \!\!\!&=&\!\!\! \{\alpha \in \Iso_C(\gJ^C) \, | \, \alpha(X \circ Y) = \alpha X \circ \alpha Y\},
\\[0.5mm]
       F_4 \!\!\!&=&\!\!\! \{\alpha \in \Iso_R(\gJ) \, |\, \alpha(X \circ Y) = \alpha X \circ \alpha Y\},
\end{eqnarray*}
respectively, and the simply connected complex Lie group ${E_6}^C$ is given by
$$
     {E_6}^C = \{\alpha \in \Iso_C(\gJ^C) \, | \, \det \, \alpha X = \det \, X \}. $$

   We define $\R$-linear transformations $\sigma$ and $\sigma'$ of $\gJ$ by
$$
   \sigma X=\sigma\begin{pmatrix} \xi_1 & x_3 & \ov{x}_2 \\
                         \ov{x}_3 & \xi_2 & x_1 \\
                         x_2 & \ov{x}_1 &\xi_3     \end{pmatrix} 
        = \begin{pmatrix} \xi_1 & -x_3 & -\ov{x}_2 \\
                         -\ov{x}_3 & \xi_2 & x_1 \\
                          -x_2 & \ov{x}_1 & \xi_3   \end{pmatrix},
 \quad \sigma'X = \begin{pmatrix} \xi_1 & x_3 & -\ov{x}_2 \\
                                  \ov{x}_3 & \xi_2 & -x_1 \\                                   
                                  -x_2 & -\ov{x}_1 & \xi_3  \end{pmatrix}, $$respectively. Then $\sigma, \sigma' \in F_4 \subset {F_4}^C$. $\sigma$ and $\sigma'$ are commutative: $\sigma\sigma' = \sigma'\sigma$.
\vspace{1mm}

Let $\mathfrak{P}^C$ be the Freudenthal $C$-vector space 
$$
     \mathfrak{P}^C = \mathfrak{J}^C \oplus \mathfrak{J}^C \oplus C \oplus C, $$
in which the Freudenthal cross operation $P \times Q, P = (X, Y, \xi, \eta), Q = (Z, W, \zeta, \omega) \in \gP^C$, is defined as follow:
$$
      P \times Q = {\varPhi}(\phi, A, B, \nu), \quad
     \left \{ \begin{array}{l}
\vspace{1mm} 
          \phi = - \dfrac{1}{2}(X \vee W + Z \vee Y) \\
\vspace{1mm}
           A = - \dfrac{1}{4}(2Y \times W - \xi Z - \zeta X) \\
\vspace{1mm}
           B =  \dfrac{1}{4}(2X \times Z - \eta W - \omega Y) \\
\vspace{1mm}
           \nu = \dfrac{1}{8}((X, W) + (Z, Y) - 3(\xi\omega + \zeta\eta)),
\end{array} \right. 
$$
\noindent where $X \vee W \in {\gge_6}^C$ is defined by 
$$
    X \vee W = [\tilde{X}, \tilde{W}] + (X \circ W - \dfrac{1}{3}(X, W)E)^{\sim}, 
$$ 
\noindent  here $\tilde{X} : \gJ^C \to \gJ^C$ is defined by $\tilde{X}Z = X \circ Z$, $Z \in \gJ^C$. 

Now, the simply connected complex Lie group ${E_7}^C$ is defined by
$$
    {E_7}^C = \{\alpha \in \Iso_C(\gP^C) \, | \, \alpha(P \times Q)\alpha^{-1} = 
\alpha P \times \alpha Q\}. $$
The Lie algebra ${\mathfrak{e}_7}^C$ of the group ${E_7}^C$ is given by
$$
    {\mathfrak{e}_7}^C = \{\varPhi(\phi, A, B, \nu) \, | \, \phi \in {\mathfrak{e}_6}^C, A, B \in \gJ^C, \nu \in C\}. $$
Naturally we have ${F_4}^C \subset {E_6}^C \subset {E_7}^C$. Finally, in a $C$-vector space ${\mathfrak{e}_8}^C$:      
$$
    {\mathfrak{e}_8}^C = {\mathfrak{e}_7}^C \oplus \gP^C \oplus \gP^C \oplus C \oplus C \oplus C, $$
we define a Lie bracket $[R_1, R_2]$ by
$$
  [(\varPhi_1, P_1, Q_1, r_1, s_1, t_1), (\varPhi_2, P_2, Q_2, r_2, s_2, t_2)] = 
  (\varPhi, P, Q, r, s, t), $$
$$
\left\{\begin{array}{l}
     \varPhi = {[}\varPhi_1, \varPhi_2] + P_1 \times Q_2 - P_2 \times Q_1
\vspace{1mm} \\
     Q = \varPhi_1P_2 - \varPhi_2P_1 + r_1P_2 - r_2P_1 + s_1Q_2 - s_2Q_1 
\vspace{1mm} \\
     P = \varPhi_1Q_2 - \varPhi_2Q_1 - r_1Q_2 + r_2Q_1 + t_1P_2 - t_2P_1
\vspace{1mm} \\
     r = - \dfrac{1}{8}\{P_1, Q_2\} + \dfrac{1}{8}\{P_2, Q_1\} + s_1t_2 - s_2t_1\vspace{1mm} \\
     s = \,\,\, \dfrac{1}{4}\{P_1, P_2\} + 2r_1s_2 - 2r_2s_1
\vspace{1mm} \\
     t = - \dfrac{1}{4}\{Q_1, Q_2\} - 2r_1t_2 + 2r_2t_1, 
\end{array} \right. $$
where $\{P, Q \} = (X, W) -(Y, Z) + \xi\omega - \eta\zeta, P = (X, Y, \xi, \eta),\,Q = (Z, W, \zeta, \omega) \in \vspace{0.5mm}\gP^C$,
then ${\mathfrak{e}_8}^C$ becomes a complex simple Lie algebra of type $E_8$ ([2]).  
\vspace{2mm}

We define a $C$-linear transformation $\wti{\lambda}$ of ${\mathfrak{e}_8}^C$
by
$$
       \wti{\lambda}(\varPhi, P, Q, r, s, t) = (\lambda\varPhi\lambda^{-1}, \lambda Q, - \lambda P, -r, -t, -s), $$
where $\lambda$ is a $C$-linear transformation of $\gP^C$ defined by
$     
    \lambda(X, Y, \xi, \eta) = (Y, - X, \eta, -\xi).\vspace{1mm} $
   
The complex conjugation in ${\mathfrak{e}_8}^C$ is denoted by $\tau$:     
$$
    \tau(\varPhi, P, Q, r, s, t) = (\tau\varPhi\tau, \tau P, \tau Q, \tau r, \tau s, \tau t),$$
where $\tau$ in the right hand side is the usual complex conjugation in the complexification.
\vspace{2mm}

Now, the connected complex Lie group ${E_8}^C$ and connected compact Lie group $E_8$ are given by 
\begin{eqnarray*}
   {E_8}^C \!\!\!&=&\!\!\! \{\alpha \in \Iso_C({\mathfrak{e}_8}^C) \,|\, \alpha[R, R'] = [\alpha R, \alpha R']\}, 
\vspace{1mm} \\    
   E_8 \!\!\!&=&\!\!\! \{\alpha \in {E_8}^C \, | \, \tau\wti{\lambda}\alpha\wti{\lambda}\tau = \alpha\} = ({E_8}^C)^{\tau\wti{\lambda}}, 
\end{eqnarray*}
respectively.
\vspace{1mm}

For $\alpha \in {E_7}^C$, the mapping $\wti{\alpha} : {\mathfrak{e}_8}^C \to {\mathfrak{e}_8}^C$ is defined by
$$
     \wti{\alpha}(\varPhi, P, Q, r, s, t) = (\alpha\varPhi\alpha^{-1}, \alpha P, \alpha Q, r, s, t),
$$
then $\wti{\alpha} \in {E_8}^C$, so $\alpha$ and $\wti{\alpha}$ will be identified. The group ${E_8}^C$ contains ${E_7}^C$ as a subgroup by
\begin{eqnarray*}
       {E_7}^C \!\!\!&=&\!\!\! \{\wti{\alpha} \in {E_8}^C \,| \, \alpha \in {E_7}^C\}
\vspace{1mm}\\
             &=&\!\!\! ({E_8}^C)_{(0,0,0,1,0,0),(0,0,0,0,1,0),(0,0,0,0,0,1)}.
\end{eqnarray*}
Similarly, $E_7 \subset E_8$. 
In particular, elements $\sigma, \sigma'$ of 
$F_4$ are also elements of $E_8 \subset {E_8}^C$. Therefore the actions of $\sigma$ and $\sigma'$ on ${\mathfrak{e}_8}^C$ are given as 
\begin{eqnarray*}
     \sigma(\varPhi, P, Q, r, s, t) \!\!\!&=&\!\!\! (\sigma\varPhi\sigma^{-1}, \sigma P, \sigma Q, r, s, t), 
\vspace{1mm}\\
     \sigma'(\varPhi, P, Q, r, s, t) \!\!\!&=&\!\!\! (\sigma'\varPhi\sigma'^{-1}, \sigma'P, \sigma'Q, r, s, t).        
\end{eqnarray*}

Hereafter, in ${\mathfrak{e}_8}^C$, we shall use the following notations.
$$
\begin{array}{c}
      \varPhi = (\varPhi, 0, 0, 0, 0, 0), \quad P^- = (0, P, 0, 0, 0, 0), \quad Q_- = (0, 0, Q, 0, 0, 0), 
\vspace{1mm}\\
    \wti{r} = (0, 0, 0, r, 0, 0), \quad s^- = (0, 0, 0, 0, s, 0), \quad 
\vspace{5mm}
t_- = (0, 0, 0, 0, 0, t). 
\end{array} $$

\section
{\bf  The group $({F_4}^C)^{\sigma, \sigma'}$}
\vspace{3mm}

We shall investigate the structure of a subgroup $({F_4}^C)^{\sigma,\sigma'}$ of the group ${F_4}^C$:
$$
       ({F_4}^C)^{\sigma,\sigma'} = \{\alpha \in {F_4}^C \,| \, \sigma\alpha = \alpha\sigma, \sigma'\alpha = 
\vspace{3mm}
\alpha\sigma'\}. $$
       
{\bf Lemma 2.1.} {\it In the Lie algebra ${\mathfrak{f}_4}^C$ of the group ${F_4}^C:$
$$
        {\mathfrak{f}_4}^C = \{\delta \in \Hom_C(\gJ^C) \, | \, \delta(X \circ Y) = \delta X \circ Y + X \circ \delta Y \}, $$
the Lie algebra $({\gf_4}^C)^{\sigma,\sigma'}$ of the group $({F_4}^C)^{\sigma,\sigma'}$
is isomorphic to $\spin(8, C)$}:
$$
    ({\gf_4}^C)^{\sigma,\sigma'} \cong \spin(8, C) = \so(8, C). $$
{\it The action of $D = (D_1, D_2, D_3) \in \spin(8, C)$ on $\gJ^C$ is given by
$$
   D\begin{pmatrix} \xi_1 & x_3 & \ov{x}_2 \\
                   \ov{x}_3 & \xi_2 & x_1 \\     
                   x_2 & \ov{x}_1 & \xi_3 \end{pmatrix} = 
      \begin{pmatrix} 0 & D_3(x_3) & \ov{D_2(x_2)} \\
                   \ov{D_3(x_3)} & 0 & D_1(x_1) \\     
                   D_2(x_2) & \ov{D_1(x_1)} & 0 \end{pmatrix}, $$
where $D_1, D_2, D_3$ are elements of $\so(8, C)$ satisfying}
$$
       D_1(x)y + xD_2(y) = \ov{D_3(\ov{xy})}, 
\vspace{3mm}
\quad x, y \in \gC^C. $$            

{\bf Theorem 2.2.} \quad \qquad $({F_4}^C)^{\sigma,\sigma'} \cong Spin(8, C)$.
\begin{proof}Let $Spin(8, C)$ be the group defined by
$$
   \{(\alpha_1, \alpha_2, \alpha_3) \in SO(8, C) \times SO(8, C) \times SO(8, C) \, |\, (\alpha_1x)(\alpha_2y) = \ov{\alpha_3(\ov{xy})}, x, y \in \gC^C \} $$
(cf.[7] Theorem 1.47). Now, we define a mapping $\varphi : Spin(8, C) \to ({F_4}^C)^{\sigma,\sigma'}$ by
$$
    \varphi(\alpha_1, \alpha_2, \alpha_3)\begin{pmatrix} \xi_1 & x_3 & \ov{x}_2 \\
                                                         \ov{x}_3 & \xi_2 & x_1 \\
                                                         x_2 & \ov{x}_1 & \xi_3 
                                                         \end{pmatrix} = 
  \begin{pmatrix} \xi_1 & \alpha_3x_3 & \ov{\alpha_2x_2} \\
                   \ov{\alpha_3x_3} & \xi_2 & \alpha_1x_1 \\     
                   \alpha_2x_2 & \ov{\alpha_1x_1} & \xi_3 \end{pmatrix}. $$ 
                        
\noindent $\varphi$ is well-defined: $\varphi(\alpha_1,\alpha_2,\alpha_3) \in ({F_4}^C)^{\sigma,\sigma'}$. Clearly $\varphi$ is a homomorphism and injective. Since $({F_4}^C)^{\sigma,\sigma'} = (({F_4}^C)^\sigma)^{\sigma'} = (Spin(9, C))^{\sigma'}$ ([6] Theorem 2.4.3) is connected and $\dim_C (({\mathfrak{f}_4}^C)^{\sigma,\sigma'}) = \dim_C (\mathfrak{spin}(8, C))$ (Lemma 2.1), $\varphi$ is onto. Thus we have the required isomorphism $({F_4}^C)^{\sigma,\sigma'} \cong Spin(8, C)$. 
\end{proof}
\vspace{3mm}

\section
{\bf The groups $({E_7}^C)^{\sigma,\sigma'}$ and $({E_7}^C)^{\sigma,\sigma',\mathfrak{so}(8,C)}$ }
\vspace{3mm}

The aim of this section is to show the connectedness of the group $({E_7}^C)^{\sigma,\sigma',\mathfrak{so}(8,C)}$\!. Now, we define subgroups $({E_7}^C)^{\sigma,\sigma'}$ and $({E_7}^C)^{\sigma,\sigma',\mathfrak{so}(8,C)}$ of  the group ${E_7}^C$ respectively by
\begin{eqnarray*}
     ({E_7}^C)^{\sigma,\sigma'} \!\!\! &=&\!\!\! \{\alpha \in {E_7}^C \, | \, \sigma\alpha = \alpha\sigma, \sigma'\alpha = \alpha\sigma' \},
\vspace{1mm}\\
     ({E_7}^C)^{\sigma,\sigma',\mathfrak{so}(8,C)} \!\!\! &=&\!\!\! \{\alpha \in ({E_7}^C)^{\sigma,\sigma'} \, | \, \varPhi_D\alpha = \alpha\varPhi_D \; \mbox{for all}\; D \in \mathfrak{so}(8, C)\},
\end{eqnarray*}
where $\varPhi_D = (D, 0, 0, 0) \in {\mathfrak{e}_7}^C, D \in \so(8, C) = ({\mathfrak{f}_4}^C)^{\sigma,\sigma'}$.
\vspace{3mm}

{\bf Lemma 3.1.} (1) {\it The Lie algebra $({\mathfrak{e}_7}^C)^{\sigma,\sigma'}$ of the group  $({E_7}^C)^{\sigma,\sigma'}$ is given by}
$$
\begin{array}{lll}
  & & ({\mathfrak{e}_7}^C)^{\sigma,\sigma'} 
\vspace{1mm}\\
  &=&\!\! \{\varPhi \in {\mathfrak{e}_7}^C \, | \, \sigma\varPhi = \varPhi\sigma, \sigma'\varPhi = \varPhi\sigma'\}
\vspace{1mm}\\
  &=&\!\! \left\{\varPhi(\phi, A, B, \nu) \in {\mathfrak{e}_7}^C \,\left| \, 
\begin{array}{l}
     \phi \in ({\mathfrak{e}_6}^C)^{\sigma,\sigma'}, A, B \in \gJ^C,
\\
     A, B \; \mbox{are diagonal forms}, \; \nu \in C 
\end{array}\right. \right\}.
\end{array}
$$
\vspace{1mm}

(2) {\it The Lie algebra $({\mathfrak{e}_7}^C)^{\sigma,\sigma',\mathfrak{so}(8,C)}$ of the group $({E_7}^C)^{\sigma,\sigma',\mathfrak{so}(8,C)}$ is given by}
$$
\begin{array}{lll}
&  &({\mathfrak{e}_7}^C)^{\sigma,\sigma',\mathfrak{so}(8,C)} 
\vspace{1mm}\\
&=&\!\!\! \{\varPhi \in ({\mathfrak{e}_7}^C)^{\sigma,\sigma'} \, | \, [\varPhi, \varPhi_D] = 0 \; \mbox{for all} \;D \in \so(8,C)\}
\vspace{1mm}\\
&=&\!\!\! \left\{\varPhi(\phi, A, B, \nu) \in ({\mathfrak{e}_7}^C)^{\sigma,\sigma'} \,\left|  \begin{array}{l}
     \phi = \wti{T}, T = \diag(\tau_1, \tau_2, \tau_3) \in \gJ^C, \tau_k \in C,
\\
     \tau_1 + \tau_2 + \tau_3  = 0, A = \diag(\alpha_1, \alpha_2, \alpha_3) \in \gJ^C, \\
\alpha_k \in C, B = \diag(\beta_1, \beta_2, \beta_3) \in \gJ^C, \beta_k \in C 
\end{array}\right.
\!\! \right \}\!\!.
\end{array}
$$

\noindent {\it In particular, we have} 
$$
     \dim_C(({\mathfrak{e}_7}^C)^{\sigma,\sigma',\mathfrak{so}(8,C)}) 
\vspace{2mm}
= 2 + 3 \times 2 + 1 = 9. $$

In the Lie algebra ${\mathfrak{e}_7}^C$, we define 
$$
      \kappa = \varPhi(-2E_1 \vee E_1, 0, 0, -1), \quad \mu = \varPhi(0, E_1, E_1, 0), $$
and we define the group $({E_7}^C)^{\kappa,\mu}$ by
$$
       ({E_7}^C)^{\kappa,\mu} = \{\alpha \in {E_7}^C \, | \, \kappa\alpha = \alpha\kappa, \mu\alpha = \alpha\mu\}. 
$$
\indent Note that if $\alpha \in {E_7}^C$ satisfies $\kappa\alpha=\alpha\kappa$, then $\alpha$ automatically satisfies $\sigma\alpha=\alpha\sigma$ because $-\sigma=\exp(\pi i \kappa), i \in C, i^2=-1$.
\vspace{2mm}
         
{\bf Lemma 3.2.} {\it For $A \in SL(2,C) = \{A \in M(2, C) \, |\, \det \, A = 1\}$, we define $C$-linear transformations $\phi_k(A), k = 1, 2, 3$ of $\gP^C$ by}
$$
\begin{array}{c}
     \phi_k(A)\Big(\begin{pmatrix} \xi_1 & x_3 & \ov{x}_2 \\
                                   \ov{x}_3 & \xi_2 & x_1 \\
                                   x_2 & \ov{x}_1 & \xi_3 \end{pmatrix},
                   \begin{pmatrix} \eta_1 & y_3 & \ov{y}_2 \\
                                   \ov{y}_3 & \eta_2 & y_1 \\
                                   y_2 & \ov{y}_1 & \eta_3 \end{pmatrix}, \xi, \eta \Big)
\vspace{1.5mm}\\
\hspace*{-5mm}
    = \Big(\begin{pmatrix} {\xi'}_1 & {x'}_3 & \ov{x'}_2 \\
                                   \ov{x'}_3 & {\xi'}_2 & {x'}_1 \\
                                   {x'}_2 & \ov{x'}_1 & {\xi'}_3 \end{pmatrix},
                   \begin{pmatrix} {\eta'}_1 & {y'}_3 & \ov{y'}_2 \\
                                   \ov{y'}_3 & {\eta'}_2 & {y'}_1 \\
                                   {y'}_2 & \ov{y'}_1 & {\eta'}_3 \end{pmatrix}, \xi', \eta' \Big),
\vspace{1.5mm}\\
\hspace*{-2mm}\begin{pmatrix}{\xi'}_k \\
                      \eta' \end{pmatrix}\! = \!A\begin{pmatrix} \xi_k \\                                                           \eta \end{pmatrix}, 
\begin{pmatrix}{\xi'} \\
                      {\eta'}_k \end{pmatrix} \!=\! A\begin{pmatrix} \xi \\                                                           \eta_k \end{pmatrix}, 
     \begin{pmatrix} {\eta'}_{k+1} \\
                      {\xi'}_{k+2} \end{pmatrix} \!=\! A\begin{pmatrix} \eta_{k+1} \\                                                           \xi_{k+2} \end{pmatrix},  
      \begin{pmatrix} {\eta'}_{k+2} \\
                      {\xi'}_{k+1} \end{pmatrix} \!=\! A\begin{pmatrix} \eta_{k+2} \\                                                          \xi_{k+1} \end{pmatrix}, 
\vspace{1.5mm}\\
\hspace*{-4mm}\begin{pmatrix}{x'}_k \\
                 {y'}_k \end{pmatrix} = {}^t\!A^{-1} \begin{pmatrix} x_k \\                                                           y_k \end{pmatrix},   
   \begin{pmatrix}{x'}_{k+1} \\
                 {y'}_{k+1} \end{pmatrix} = \begin{pmatrix} x_{k+1} \\                                                           y_{k+1} \end{pmatrix},   
  \begin{pmatrix}{x'}_{k+2} \\
                 {y'}_{k+2} \end{pmatrix} = \begin{pmatrix} x_{k+2} \\                                                           y_{k+2} \end{pmatrix}   
\end{array}$$           
({\it where indices are considered as} mod 3). {\it Then $\phi_k(A) \!\in \! ({E_7}^C)^{\sigma,\sigma',\mathfrak{so}(8,C)}\! \subset \!({E_7}^C)^{\sigma,\sigma'}$\!. Moreover, $\phi_k(A) \in (({E_7}^C)^{\kappa,\mu})^{\sigma',\mathfrak{so}(8,C)} \subset (({E_7}^C)^{\kappa,\mu})^{\sigma'}$ for $k = 2, 3$.}              
\vspace{5mm}

{\bf Proposition 3.3.} $({E_7}^C)^\sigma \cong (SL(2, C) \times Spin(12, C))/\Z_2, \Z_2 = \{(E, 1), (-E, $ $-\sigma)\}$.
\begin{proof} Let $Spin(12, C) = ({E_7}^C)^{\kappa,\mu}$ ([6] Proposition 4.6.10) $\subset ({E_7}^C)^\sigma$.
We define a mapping $\varphi_1 : SL(2, C) \times Spin(12, C) \to ({E_7}^C)^\sigma$ by
$$
     \varphi_1(A_1, \delta) = \phi_1(A_1)\delta. $$
Then we have this proposition (see [6] Theorem 4.6.13 for details). 
\end{proof}
\vspace{3mm}

{\bf Lemma 3.4.} {\it The Lie algebra $(({\mathfrak{e}_7}^C)^{\kappa,\mu})^{\sigma'}$ of the group $(({E_7}^C)^{\kappa,\mu})^{\sigma'}$ is given by} 
$$
\begin{array}{lll}
  & &\!\!\!(({\mathfrak{e}_7}^C)^{\kappa,\mu})^{\sigma'} 
\vspace{1.5mm}\\
  \!\!\!&=&\!\!\! \{\varPhi \in {\mathfrak{e}_7}^C \, | \, \kappa\varPhi = \varPhi\kappa, \mu\varPhi = \varPhi\mu, \sigma'\varPhi = \varPhi\sigma'\}
\vspace{1.5mm}\\
  \!\!\!&=&\!\!\! \left\{\varPhi(\phi, A, B, \mu) \in {\mathfrak{e}_7}^C \, \left| \, 
          \begin{array}{l} \phi \in ({\mathfrak{e}_6}^C)^{\sigma,\sigma'}, A, B \in \gJ^C, \sigma A = \sigma'A = A, 
\\[0.5mm]
      (E_1, A) = 0, \sigma B = \sigma'B = B, (E_1, B) = 0, 
\\[0.5mm]
       \nu = -\dfrac{3}{2}(\phi E_1, E_1) \end{array} \right. \right\}.
\end{array}
$$    
       
 \noindent {\it In particular, we have}
$$
      \dim_C((({\mathfrak{e}_7}^C)^{\kappa,\mu})^{\sigma'}) 
\vspace{3mm}
= 30 + 2 \times 2 = 34. $$ 

{\bf Proposition 3.5.} $(({E_7}^C)^{\kappa,\mu})^{\sigma'} \cong (SL(2, C) \times SL(2, C) \times Spin(8, C))/\Z_2, \Z_2 $  $= \{(E, E, 1), (-E, -E, \sigma)\}$. 
\begin{proof}Let $Spin(8, C) \!=\! ({F_4}^C)^{\sigma,\sigma'}\! = \!(Spin(9, C))^{\sigma'} \! \subset \! (Spin(12, C))^{\sigma'}\! =\! (({E_7}^C)^{\kappa,\mu})^{\sigma'}$\!\!. Now, we define a mapping $\varphi_2 : SL(2, C) \times SL(2, C) \times Spin(8, C) \to (({E_7}^C)^{\kappa,\mu})^{\sigma'}$ by
$$
    \varphi_2(A_2, A_3, \beta) = \phi_2(A_2)\phi_3(A_3)\beta. 
$$
$\varphi_2$ is well-defined: $\varphi_2(A_2, A_3, \beta) \in (({E_7}^C)^{\kappa,\mu})^{\sigma'}$ (Lemma 3.2). Since $\phi_2(A_2),$ $ \phi_3(A_3)$ and $\beta$ commute with each other, $\varphi_2$ is a homomorphism. $\Ker \, \varphi_2 \,= \,\{(E, E, 1), $ $(-E,-E, \sigma)\} \,=\, \Z_2$. Since $(({E_7}^C)^{\kappa,\mu})^{\sigma'}\, = \,(Spin(12, C))^{\sigma'}$ is connected and $\dim_C((({\mathfrak{e}_7}^C)^{\kappa,\mu})^{\sigma'})$ $ = 34 $(Lemma 3.4) $= 3 + 3+ 28 = \dim_C(\mathfrak{sl}(2, C) \oplus \mathfrak{sl}(2, C) \oplus \mathfrak{so}(8, C))$, $\varphi_2$ is onto. Thus we have the required isomorphism $(({E_7}^C)^{\kappa,\mu})^{\sigma'} \cong (SL(2, C) \times SL(2, C) \times Spin(8, C))/\Z_2$. 
\end{proof}    
\vspace{3mm}

{\bf Theorem 3.6.} $({E_7}^C)^{\sigma,\sigma'} \cong (SL(2, C) \times SL(2, C) \times SL(2, C) \times Spin(8, C))$ $/(\Z_2 \times \Z_2), \Z_2 \times \Z_2 = \{(E, E, E, 1), (E, -E , -E, \sigma)\} \times \{(E, E, E, 1), (-E, -E,$ $ E, \sigma')\}$.
\begin{proof} Let $Spin(8, C) = ({F_4}^C)^{\sigma,\sigma'}$ (Theorem 2.2) $\subset ({E_7}^C)^{\sigma,\sigma'}$. We define a mapping $\varphi : SL(2, C) \times SL(2, C) \times SL(2, C) \times Spin(8, C) \to ({E_7}^C)^{\sigma,\sigma'}$ by
$$
       \varphi(A_1, A_2, A_3, \beta) = \phi_1(A_1)\phi_2(A_2)\phi_3(A_3)\beta. $$
$\varphi$ is well-defined: $\varphi(A_1, A_2, A_3, \beta) \in ({E_7}^C)^{\sigma,\sigma'}$ (Lemma 3.2). Since $\phi_1(A_1),$ $ \phi_2(A_2),$ $ \phi_3(A_3)$ and $\beta$ commute with each other, $\varphi$ is a homomorphism. We shall show that $\varphi$ is 
onto. For $\alpha \in ({E_7}^C)^{\sigma,\sigma'} \subset ({E_7}^C)^\sigma$, there exist $A_1 \in SL(2, C)$ and $\delta \in Spin(12, C)$ such that $\alpha = \phi_1(A_1)\delta$ (Proposition 3.3). From the condition $\sigma'\alpha\sigma' = \alpha$, that is, $\sigma'\varphi_1(A_1, \delta)\sigma' =\varphi_1(A_1, \delta)$, we have $\varphi_1(A_1, \sigma' \delta \sigma')=\varphi_1(A_1, \delta)$. Hence
$$
     \left\{\begin{array}{l}
                  A_1 = A_1
\vspace{1mm}\\
                  \sigma'\delta\sigma' = \delta, \end{array} \right.  \quad \mbox{or} \quad 
     \left\{\begin{array}{l} A_1 = -A_1
\vspace{1mm}\\
                  \sigma'\delta\sigma' = - \sigma\delta. \end{array} \right. $$
The latter case is impossible because $A_1 = -A_1$ is false. In the first case, from $\sigma'\delta\sigma' = \delta$, we have $\delta \in (Spin(12, C))^{\sigma'} = (({E_7}^C)^{\kappa,\mu})^{\sigma'}$. Hence there exist $A_2, A_3 \in SL(2, C)$ and $\beta \in Spin(8, C)$ such that $\delta = \phi_2(A_2)\phi_3(A_3)\beta$ (Proposition 3.5). Then
$$
    \alpha = \phi_1(A_1)\delta = \phi_1(A_1)\phi_2(A_2)\phi_3(A_3)\beta = \varphi(A_1, A_2, A_3, \beta). $$
Therefore $\varphi$ is onto. It is not difficult to see that
\begin{eqnarray*}
    \Ker \, \varphi &=& \{(E, E, E, 1),(E, -E, -E, \sigma), (-E, E, -E, \sigma\sigma'),(-E, -E, E, \sigma')\}
\vspace{1mm}\\
      &=& \{(E, E, E, 1),(E, -E, -E, \sigma)\} \times \{(E, E, E, 1),(-E, -E, E, \sigma')\}
\vspace{1mm}\\  
      &=& \Z_2 \times \Z_2.
\end{eqnarray*}      
Thus we have the required isomorphism $({E_7}^C)^{\sigma,\sigma'}\!\! \cong (SL(2, C) \times SL(2,C) \times SL(2, C) $ $\times Spin(8, C))$ $/(\Z_2 \times \Z_2)$.
\end{proof}
\vspace{3mm}

Now, we determine the structure of the subgroup $({E_7}^C)^{\sigma,\sigma',\mathfrak{so}(8,C)}$ of ${E_7}^C$.
\vspace{3mm}

{\bf Theorem 3.7.} \quad $({E_7}^C)^{\sigma,\sigma',\mathfrak{so}(8,C)} \cong SL(2, C) \times SL(2, C) \times SL(2, C)$.
\vspace{1mm}

\noindent {\it In particular, $({E_7}^C)^{\sigma,\sigma',\mathfrak{so}(8,C)}$ is connected.}

\begin{proof}We define a mapping $\varphi : SL(2, C) \times SL(2, C) \times SL(2, C) \to ({E_7}^C)^{\sigma,\sigma',\mathfrak{so}(8,C)}$ by
$$
       \varphi(A_1, A_2, A_3) = \phi_1(A_1)\phi_2(A_2)\phi_3(A_3). $$
($\varphi$ is the restricted mapping of $\varphi$ of Theorem 3.6). $\varphi$ is 
well-defined: $\varphi(A_1, A_2,$ $ A_3) \in ({E_7}^C)^{\sigma,\sigma',\mathfrak{so}(8,C)}$ (Lemma 3.2).
Since $\dim_C(\sl(2, C) \oplus \sl(2, C) \oplus \sl(2, C)) = 3 + 3 + 3 = 9 = \dim_C(({\mathfrak{e}_7}^C)^{\sigma,\sigma',\mathfrak{so}(8,C)})$ (Lemma 3.1), 
$\Ker \, \varphi$ is discrete. Hence $\Ker \, \varphi$ is contained in the center 
$z(SL(2, C) \times SL(2, C) \times SL(2, C)) = \{(\pm E, \pm E,$ $ \pm E)\}$. 
However, $\varphi$ maps these elements to $\pm 1, \pm \sigma, \pm \sigma', \pm \sigma\sigma'$. Hence $\Ker \, \varphi = \{(E, E, E)\}$, that is, $\varphi$ is injective. Finally, we shall show that $\varphi$
is onto. For $\alpha \in ({E_7}^C)^{\sigma,\sigma',\mathfrak{so}(8,C)} \subset ({E_7}^C)^{\sigma,\sigma'}$, there exist $A_1, A_2, A_3 \in SL(2, C)$ and 
$\beta \in Spin(8, C)$ such that $\alpha = \varphi(A_1, A_2, A_3, \beta)$ 
(Theorem 3.6). From $\varPhi_D\alpha = \alpha\varPhi_D$, we have $\varPhi_D\beta = \beta\varPhi_D$ for all $D \in \so(8,C)$, that is, $\beta$ is contained in the center $z(({F_4}^C)^{\sigma,\sigma'}) = z(Spin(8, C)) = \{1, \sigma,\sigma', \sigma\sigma'\}$. However, $\sigma = \phi_1(E)\phi_2(-E)\phi_3(-E), \sigma' = \phi_1(-E)$ $\phi_2(-E)\phi_3(E)$, so $1, \sigma, \sigma', \sigma\sigma' \in \phi_1(SL(2, C))$ $\phi_2(SL(2, C))\phi_3(SL(2, C))$. Hence we see that $\varphi$ is onto. Thus we have the required isomorphism 
$({E_7}^C)^{\sigma,\sigma',\mathfrak{so}(8,C)} \cong SL(2, C) \times SL(2, C) $  $ \times SL(2, C)$. 
\end{proof} 
\vspace{5mm}

\section
{\bf  Connectedness of the group $({E_8}^C)^{\sigma,\sigma',\mathfrak{so}(8,C)}$ }
\vspace{3mm}


We define a subgroup  $({E_8}^C)^{\sigma,\sigma',\mathfrak{so}(8,C)}$ of the group $({E_8}^C)^{\sigma,\sigma'}$ by
$$
     ({E_8}^C)^{\sigma,\sigma',\mathfrak{so}(8,C)} = \left\{\alpha
\in {E_8}^C \, \left| \, \begin{array}{l} \sigma\alpha = \alpha\sigma, \,\sigma'\alpha = \alpha\sigma', \\
 \varTheta(R_D)\alpha = \alpha\varTheta(R_D) \; \mbox{for all}\; D \in \so(8, C) \end{array} \right. \right\}, $$
where $R_D =  (\varPhi_D, 0, 0, 0, 0, 0) \in {\mathfrak{e}_8}^C$ and $\varTheta(R_D)$ means $\ad(R_D)$. Hereafter for $R \in {\mathfrak{e}_8}^C$, we denote $\ad(R)$ by ${\varTheta}(R)$.
\vspace{2mm}

To prove the connectedness of the group $({E_8}^C)^{\sigma,\sigma',\mathfrak{so}(8,C)}$, we use the method used in [5]. However, we write this method in detail again.
Firstly, we consider a subgroup $(({E_8}^C)^{\sigma,\sigma',\mathfrak{so}(8,C)})_{1_-}$ of $({E_8}^C)^{\sigma,\sigma',\mathfrak{so}(8,C)}$:
$$
   (({E_8}^C)^{\sigma,\sigma',\mathfrak{so}(8,C)})_{1_-} = \{\alpha \in ({E_8}^C)^{\sigma,\sigma',\mathfrak{so}(8,C)} \, | \, \alpha1_- = 1_- \}. 
$$
\vspace{1mm}

{\bf Lemma 4.1.} (1) {\it The Lie algebra $(({\mathfrak{e}_8}^C)^{\sigma,\sigma',\mathfrak{so}(8,C)})_{1_-}$ of the group \\ $(({E_8}^C)^{\sigma,\sigma',\mathfrak{so}(8,C)})_{1_-}$ is given by}
$$
\begin{array}{lll}
&  &(({\mathfrak{e}_8}^C)^{\sigma,\sigma',\mathfrak{so}(8,C)})_{1_-}
\vspace{2mm}\\
&=&\left\{ R \in {\mathfrak{e}_8}^C \, \Bigg| \, 
\begin{array}{l}
    \sigma R=R, \, {\sigma'}R=R,
\\[1mm]
       [R_D, R] =0 \,\, \text{for all} \,\, D \in \mathfrak{so}(8,C),
       [R, 1_-] =0 
\end{array}   \right\} 
\vspace{2mm}\\
&=& \left\{(\varPhi, 0, Q, 0, 0, t) \in {\mathfrak{e}_8}^C \, \left| \,
\begin{array}{l}
     \varPhi \in ({\mathfrak{e}_7}^C)^{\sigma,\sigma',\mathfrak{so}(8,C)}, Q = (Z, W, \zeta, \omega), \\
     Z, W \; \mbox{are diagonal forms}, \; \zeta, \omega,t \in C
\end{array} \right. \right\}. 
\end{array}
$$
\vspace{1mm}

\noindent {\it In particular, }
$$
      \dim_C((({\mathfrak{e}_8}^C)^{\sigma,\sigma',\mathfrak{so}(8,C)})_{1_-}) = 
\vspace{3mm}
9 + 8 + 1 = 18. 
$$
(2) {\it The Lie algebra $({\mathfrak{e}_8}^C)^{\sigma,\sigma',\mathfrak{so}(8,C)}$ of the group $({E_8}^C)^{\sigma,\sigma',\mathfrak{so}(8,C)}$ is given by}
$$
\begin{array}{lll}
& & ({\mathfrak{e}_8}^C)^{\sigma, \sigma', \mathfrak{so}(8,C)}
\vspace{2mm}\\
&=&\!\! \biggl\{ R \in {\mathfrak{e}_8}^C \,\biggm|\, 
\begin{array}{l}
     \sigma R= R, \sigma' R= R,        
  \\[1mm]
    [R_D, R] =0 \,\, \text{for all} \,\, D \in \so(8,C)                                           
\end{array} \biggr \}
\vspace{2mm}\\
 & = &\!\! \left\{ (\varPhi,P, Q, r, s, t) \in {\mathfrak{e}_8}^C \,\left|\,\begin{array}{l} 
           \varPhi \in ({\mathfrak{e}_7}^C)^{\sigma, \sigma',\mathfrak{so}(8,C)},  
        \vspace{1mm}\\
           P=(X, Y,\xi, \eta), Q=(Z, W, \zeta, \omega),
        \vspace{1mm}\\
          X, Y, Z, W \text{\,are diagonal forms},
        \vspace{1mm}\\
          \xi,\eta, \zeta, \omega, r, u, v \in C 
                                              \end{array}\right. \right \}. 
\end{array}
$$
In the following proposition, we denote by $(\gP^C)_d$ the subspace of $\gP^C$:
$$
   (\gP^C)_d = \{(X, Y, \xi, \eta) \in \gP^C \, | \, X, Y \; \mbox{are diagonal forms}, 
\vspace{3mm}
\xi, \eta \in C \}. $$ 

{\bf Proposition 4.2.} {\it The group $(({E_8}^C)^{\sigma,\sigma',\mathfrak{so}(8,C)})_{1_-}$ is a semi-direct product of groups $\exp(\varTheta(((\gP^C)_d) \oplus C_- ))$ and }$({E_7}^C)^{\sigma,\sigma',\mathfrak{so}(8,C)}$:
$$
  (({E_8}^C)^{\sigma,\sigma',\mathfrak{so}(8,C)})_{1_-} = \exp(\varTheta(((\gP^C)_d) \oplus C_- ))\cdot({E_7}^C)^{\sigma,\sigma',\mathfrak{so}(8,C)}. $$

\noindent {\it In particular, $(({E_8}^C)^{\sigma,\sigma',\mathfrak{so}(8,C)})_{1_-}$ is connected}.
\begin{proof}Let $((\gP^C)_d)_- \oplus C_-) = \{(0, 0, Q, 0, 0, t) \, | \, Q \in (\gP^C)_d, t \in C\}$ be the subalgebra of 
$(({\mathfrak{e}_8}^C)^{\sigma,\sigma',\mathfrak{so}(8,C)})_{1_-}$ (Lemma 4.1.(1)). From $[Q_-, t_-] = 0$, $\varTheta(Q_-)$ commutes with $\varTheta(t_-)$. Hence we have $\exp(\varTheta(Q_- + t_-)) = \exp(\varTheta(Q_-))\exp(\varTheta(t_-))$.
Therefore $\exp(\varTheta(((\gP^C)_d)_- \oplus C_-))$ is a subgroup of $(({E_8}^C)^{\sigma,\sigma',\mathfrak{so}(8,C)})_{1-}$. Now, let $\alpha \in (({E_8}^C)^{\sigma,\sigma',\mathfrak{so}(8,C)})_{1-}$ and set
$$
    \alpha\wti{1} = (\varPhi, P, Q, r, s, t), \quad \alpha1^- = (\varPhi_1, P_1, Q_1, r_1, s_1, t_1). $$
Then from the relation $[\alpha\wti{1}, 1_-] = \alpha[\wti{1}, 1_-] = -2\alpha1_- = -21_-, [\alpha1^-, 1_-] = \alpha[1^-, 1_-] $ $= \alpha\wti{1}$, we have
$$
    P = 0, \; s = 0, \; r = 1, \; \varPhi = 0, \; P_1 = - Q, \; s_1 = 1, \; r_1 = -\dfrac{t}{2}. $$
Moreover, from $[\alpha\wti{1}, \alpha1^-] = \alpha[\wti{1}, 1^-] = 2\alpha1^-$, we have
$$
    \varPhi_1 = \dfrac{1}{2}Q \times Q, \; Q_1 = - \dfrac{t}{2}Q - \dfrac{1}{3}\varPhi_1Q, \; t_1 = -\dfrac{t^2}{4} - \dfrac{1}{16}\{Q, Q_1\}. $$
So, $\alpha$ is of the form
$$
     \alpha = \begin{pmatrix} * & * & * & 0 & \dfrac{1}{2}Q \times Q & 0 \\
                              * & * & * & 0 & -Q & 0 \\
                   * & * & * & Q & -\dfrac{t}{2}Q -\dfrac{1}{6}(Q \times Q)Q & 0 \\
                              * & * & * & 1 & -\dfrac{t}{2} & 0 \\  
                              * & * & * & 0 & 1 & 0 \\
          * & * & * & t & -\dfrac{t^2}{4} + \dfrac{1}{96}\{Q, (Q \times Q)Q\} & 1
                              \end{pmatrix}. $$
On the other hand, we have
\begin{eqnarray*}
     \delta1^- &=& \exp\Big(\varTheta(\Big(\dfrac{t}{2}\Big)_-)\Big)\exp(\varTheta(Q_-))1^-
\vspace{1mm}\\
        &=& \begin{pmatrix} \dfrac{1}{2}Q \times Q
\vspace{1mm}\\
                 - Q 
\vspace{1mm}\\
            - \dfrac{t}{2}Q - \dfrac{1}{6}(Q \times Q)Q 
\vspace{1mm}\\
           -\dfrac{t}{2} 
\vspace{1mm}\\
           1
\vspace{1mm}\\
       -\dfrac{t^2}{4} + \dfrac{1}{96}\{Q, (Q \times Q)Q\} \end{pmatrix} = \alpha1^-,
\end{eqnarray*}           
and also we get
$$
        \delta\wti{1} = \alpha\wti{1}, \;\; \delta1_- = \alpha1_- $$
easily. Hence we see that $\delta^{-1}\alpha \in (({E_8}^C)^{\sigma,\sigma', \mathfrak{so}(8,C)})_{\wti{1},1^-,1_-} = ({E_7}^C)^{\sigma,\sigma',\mathfrak{so}(8,C)}$. Therefore we have
$$
     (({E_8}^C)^{\sigma,\sigma',\mathfrak{so}(8,C)})_{1_-} = \exp(\varTheta(((\gP^C)_d)_- \oplus C_-))({E_7}^C)^{\sigma,\sigma',\mathfrak{so}(8,C)}.  $$      
Furthermore, for $\beta \in ({E_7}^C)^{\sigma,\sigma',\mathfrak{so}(8,C)}$, it is easy to see that
$$
   \beta(\exp(\varTheta(Q_-)))\beta^{-1} = \exp(\varTheta(\beta Q_-)),\quad \beta((\exp(\varTheta(t_-)))\beta^{-1} = \exp(\varTheta(t_-). $$
This shows that $\exp(\varTheta(((\gP^C)_d)_-) \oplus C_-)) = \exp(\varTheta(((\gP^C)_d)_-)\exp(\varTheta(C_-))$ is a normal subgroup of $(({E_8}^C)^{\sigma,\sigma',\mathfrak{so}(8,C)})_{1_-}$. Moreover, we have a split exact sequence
$$
   1 \to \exp(\varTheta(((\gP^C)_d)_- \oplus C_-)) \to (({E_8}^C)^{\sigma,\sigma',\mathfrak{so}(8,C)})_{1_-} \to ({E_7}^C)^{\sigma,\sigma',\mathfrak{so}(8,C)} \to 1. $$
Hence $(({E_8}^C)^{\sigma,\sigma',\mathfrak{so}(8,C)})_{1_-}$ is a semi-direct product of $\exp(\varTheta(((\gP^C)_d)_- \oplus C_-))$ and $({E_7}^C)^{\sigma,\sigma',\mathfrak{so}(8,C)}$:
$$
    (({E_8}^C)^{\sigma,\sigma',\mathfrak{so}(8,C)})_{1_-} = 
\exp(\varTheta(((\gP^C)_d)_- \oplus C_-))\cdot({E_7}^C)^{\sigma,\sigma',\mathfrak{so}(8,C)}. $$    
$\exp(\varTheta(((\gP^C)_d)_- \oplus C_-))$ is connected and $({E_7}^C)^{\sigma,\sigma',\mathfrak{so}(8,C)}$ is connected (Theorem 3.7), hence $(({E_8}^C)^{\sigma,\sigma',\mathfrak{so}(8,C)})_{1_-}$ is also connected. 
\end{proof}
\vspace{3mm}

For $R \in {\mathfrak{e}_8}^C$, we define a $C$-linear mapping $R \times R : {\mathfrak{e}_8}^C \to {\mathfrak{e}_8}^C$ by
$$
     (R \times R)R_1 = [R, [R, R_1]\,] + \dfrac{1}{30}B_8(R, R_1)R, \quad R_1 \in {\mathfrak{e}_8}^C $$
(where $B_8$ is the Killing form of the Lie algebra ${\mathfrak{e}_8}^C$) and a space $\gW^C$ by
$$
     \gW^C = \{R \in {\mathfrak{e}_8}^C \, | \, R \times R = 0, R \not= 0\}. $$
Moreover, we define a subspace $(\gW^C)_{\sigma,\sigma',\mathfrak{so}(8,C)}$ of $\gW^C$ by
$$
      (\gW^C)_{\sigma,\sigma',\mathfrak{so}(8,C)} = \{R \in \gW^C \, | \, \sigma R = R, \sigma'R = R, [R_D, R] = 0 \; \mbox{for all}\; D \in 
\vspace{3mm}
\mathfrak{so}(8,C)\}. $$

{\bf Lemma 4.3.} {\it For $R = (\varPhi, P, Q, r, s, t) \in {\mathfrak{e}_8}^C$ satisfying $\sigma R = R, \sigma'R = R$ and $[R_D, R] = 0$ for all $D \in \mathfrak{so}(8,C), R \not= 0$, $R$ belongs to $(\gW^C)_{\sigma,\sigma',\mathfrak{so}(8,C)}$ if and only if $R$ satisfies the following conditions}.
\vspace{1mm}

(1) $2s\varPhi - P \times P = 0$ \quad (2) $2t\varPhi + Q \times Q = 0$ \quad (3) $2r\varPhi + P \times Q = 0$  
\vspace{1.5mm}

(4) $\varPhi P - 3rP - 3sQ = 0 $ \quad (5) $\varPhi Q + 3rQ - 3tP = 0 $
\vspace{1.5mm}

(6) $\{P, Q\} - 16(st + r^2) = 0$
\vspace{1.5mm}

(7) $2(\varPhi P \times Q_1 + 2P \times \varPhi Q_1 - rP \times Q_1 - sQ \times Q_1) - \{P, Q_1\}\varPhi = 0$
\vspace{1.5mm}

(8) $2(\varPhi Q \times P_1 + 2Q \times \varPhi P_1 + rQ \times P_1 - tP \times P_1) \!- \{Q, P_1\}\varPhi = 0$
\vspace{1.5mm}

(9) $8((P \times Q_1)Q - stQ_1 - r^2Q_1 - \varPhi^2Q_1 + 2r\varPhi Q_1) + 5\{P, Q_1\}Q 
-2\{Q, Q_1\}P = 0$
\vspace{-2.6mm}

\hspace*{-1.7mm}(10) $8((Q \times P_1)P + stP_1 + r^2P_1 + \varPhi^2P_1 + 2r\varPhi P_1) \,+ \,5\{Q, P_1\}P 
-2\{P, Q_1\}Q= 0$
\vspace{1.5mm}

\hspace*{-1.7mm}(11) $18(\ad\,\varPhi)^2\varPhi_1 + Q \times \varPhi_1P - P \times \varPhi_1Q) +   B_7(\varPhi, \varPhi_1)\varPhi = 0$
\vspace{1.5mm}

\hspace*{-1.7mm}(12) $18(\varPhi_1\varPhi P -2\varPhi\varPhi_1P - r\varPhi_1P - s\varPhi_1Q) + B_7(\varPhi, \varPhi_1)P = 0$
\vspace{1.5mm}

\hspace*{-1.7mm}(13) $18(\varPhi_1\varPhi Q -2\varPhi\varPhi_1Q + r\varPhi_1Q - t\varPhi_1P) + B_7(\varPhi, \varPhi_1)Q = 0,$

\vspace{1.5mm}
\noindent ({\it where $B_7$ is the Killing form of the Lie algebra ${\mathfrak{e}_7}^C$}) {\it for all $\varPhi_1 \in {\mathfrak{e}_7}^C, P_1, Q_1 \in \gP^C$.}
\vspace{2mm}

We denote the connected component of $({E_8}^C)^{\sigma,\sigma',\mathfrak{so}(8,C)}$ containing the unit element by $(({E_8}^C)^{\sigma,\sigma',\mathfrak{so}(8,C)})_0$. 
\vspace{3mm}

{\bf Proposition 4.4.} {\it The group $(({E_8}^C)^{\sigma,\sigma',\mathfrak{so}(8,C)})_0$ acts on $(\gW^C)_{\sigma,\sigma',\mathfrak{so}(8,C)}$ transitively.}
\vspace{-4mm}

\begin{proof} Since $\alpha \in ({E_8}^C)^{\sigma,\sigma',\mathfrak{so}(8,C)}$ leaves invariant the Killing form $B_8$ of ${\mathfrak{e}_8}^C\!: \!B_8(\alpha R_1,$ $\alpha R_2) = B_8(R_1, R_2)$, $R_k \in {\mathfrak{e}_8}^C, k = 1, 2$, we have $\alpha R \in (\gW^C)_{\sigma,\sigma',\mathfrak{so}(8,C)}$, for $R \in (\gW^C)_{\sigma,\sigma',\mathfrak{so}(8,C)}$. 
In fact,
\begin{eqnarray*}
   (\alpha R \times \alpha R)R_1 \!\!\!&=&\!\!\! [\alpha R, [\alpha R, \alpha R_1]\,] + 
(1/30)B_8(\alpha R, R_1)\alpha R
\\[0.5mm]
      \!\!\!&=&\!\!\! \alpha[\,[R, [R, \alpha^{-1}R_1]\,] + 
(1/30)B_8(R, \alpha^{-1}R_1)\alpha R
\\[0.5mm]
      \!\!\!&=&\!\!\! \alpha((R \times R)\alpha^{-1}R_1      
\\[0.5mm]
      &=& 0,
\\[0.5mm]
[R_D, \alpha R] \!\!\!&=&\!\!\! \alpha[\alpha^{-1}R_D, R] = \alpha[R_D, R] = 0.
\end{eqnarray*} 
\vspace{-4mm}
     
This shows that $({E_8}^C)^{\sigma,\sigma',\mathfrak{so}(8,C)}$ acts on $(\gW^C)_{\sigma,\sigma',\mathfrak{so}(8,C)}$. We will show that this action is transitive. 
Firstly for    
$R_1 \in {\mathfrak{e}_8}^C$, since
\begin{eqnarray*}
     (1_- \times 1_-)R_1 \!\!\!&=&\!\!\! [1_-,[1_-, (\varPhi_1, P_1, Q_1, r_1, s_1, t_1)] \, ] + 
(1/30)B_8(1_-, R_1)1_-
\\[0.5mm]
 \!\!\!&=&\!\!\! [1_-, (0, 0, P_1, -s_1, 0, 2r_1)] + 2s_11_- 
\\[0.5mm]
 \!\!\!&=&\!\!\! (0, 0, 0, 0, -2s_1) + 2s_11_- 
\\[0.5mm]
 \!\!\!&=&\!\!\! 0,   
\\[0.5mm]
 [R_D, 1_-] \!\!\!&=&\!\!\! 0,
\end{eqnarray*}


\noindent we have $1_- \in (\gW^C)_{\sigma,\sigma',\mathfrak{so}(8,C)}$. In order to prove the transitivity, it is sufficient to show that any element $R \in (\gW^C)_{\sigma,\sigma',\mathfrak{so}(8,C)}$ can be transformed to $1_- \in (\gW^C)_{\sigma,\sigma',\mathfrak{so}(8,C)}$ by a certain $\alpha \in ({E_8}^C)^{\sigma,\sigma',\mathfrak{so}(8,C)}$.
\vspace{1mm}

Case (1). $R = (\varPhi, P, Q, r, s,t), t \not= 0$. From (2),(5),(6) of Lemma 4.3, we have
$$
   \varPhi = -\dfrac{1}{2t}Q \times Q, \; P = \dfrac{r}{t}Q - \dfrac{1}{6t^2}(Q \times Q)Q, \; s = -\dfrac{r^2}{t} + \dfrac{1}{96t^3}\{Q, (Q \times Q)Q\}. 
$$
Now, for $\varTheta = \varTheta(0, P_1, 0, r_1, s_1, 0) \in \varTheta(({\gge_8}^C)^{\sigma,\sigma',\mathfrak{so}(8,C)})$ (Lemma 4.1.(2)), we compute $\varTheta^n1_-$,
\vspace{1mm}

{\small $$
\begin{array}{l}
   \varTheta^n1_-
\vspace{1mm}\\
  = \begin{pmatrix} ((-2)^{n-1} + (-1)^n){r_1}^{n-2}P_1 \times P_1 
\vspace{1mm}\\
 \Big((-2)^{n-1} - \dfrac{1 + (-1)^{n-1}}{2}\Big){r_1}^{n-2}s_1P_1 + \Big(\dfrac{1 - (-2)^n}{6} + \dfrac{(-1)^n}{2}\Big){r_1}^{n-3}(P_1 \times P_1)P_1
\vspace{1mm}\\
      ((-2)^n + (-1)^{n-1}){r_1}^{n-1}P_1
\vspace{1mm}\\
      (-2)^{n-1}{r_1}^{n-1}s_1
\vspace{1mm}\\
 -((-2)^{n-2} + 2^{n-2}){r_1}^{n-2}{s_1}^2 + \dfrac{2^{n-2} + (-2)^{n-2} - (-1)^{n-1}}{24}{r_1}^{n-4}\{P_1,(P_1 \times P_1)P_1\}
\vspace{1mm}\\
     (-2)^n{r_1}^n \end{pmatrix}. 
\end{array} $$}  
Hence by simple computing, we have
{\small 
$$
\begin{array}{l}
    \exp(\varPhi(0, P_1, 0, r_1, s_1, 0))1_- = (\exp\varTheta)1_- = \Big(\dsum_{n= 0}^\infty\dfrac{1}{n!}\varTheta^n \Big)1_-
\vspace{1mm}\\
= \begin{pmatrix}   -\dfrac{1}{2{r_1}^2}(e^{-2r_1} -2e^{-r_1} + 1)P_1 \times P_1
\vspace{1mm}\\
  \dfrac{s_1}{2{r_1}^2}(-e^{-2r_1} - e^{r_1} + e^{-r_1} + 1)P_1 + \dfrac{1}{6{r_1}^3}(-e^{-2r_1} + e^{r_1} + 3e^{-r_1} - 3)(P_1 \times P_1)P_1
\vspace{1mm}\\
   \dfrac{1}{r_1}(e^{-2r_1} - e^{-r_1})P_1
\vspace{1mm}\\  
   \dfrac{s_1}{2r_1}(1 - e^{-2r_1})
\vspace{1mm}\\
   -\dfrac{{s_1}^2}{4{r_1}^2}(e^{-2r_1} + e^{2r_1} -2) + \dfrac{1}{96{r_1}^4}(e^{2r_1} + e^{-2r_1} - 4e^{r_1} - 4e^{-r_1} + 6)\{P_1, (P_1 \times P_1)P_1\}
\vspace{1mm}\\
e^{-2r_1} \end{pmatrix}. 
\end{array} $$}   
\Big(if $r_1 = 0, \dfrac{f(r_1)}{{r_1}^k}$ means $\dlim_{r_1 \to 0}\dfrac{f(r_1)}{{r_1}^k}$ \Big). Here we set
$$
   Q = \dfrac{1}{r_1}(e^{-2r_1} - e^{-r_1})P_1, \;\; r = \dfrac{s_1}{2r_1}(1 - e^{-2r_1}), \;\; t = e^{-2r_1}. $$
Then we have
$$
    (\exp\,\varTheta)1_- = \begin{pmatrix} -\dfrac{1}{2t}Q \times Q
\vspace{1mm}\\
                                \dfrac{r}{t}Q - \dfrac{1}{6t^2}(Q \times Q)Q \vspace{1mm}\\
                              Q \vspace{1mm}\\
                              r \vspace{1mm}\\
         -\dfrac{r^2}{t} + \dfrac{1}{96 t^3}\{Q, (Q \times Q)Q\}\vspace{1mm}\\
                              t 
\end{pmatrix} = \begin{pmatrix} \varPhi \vspace{1mm}\\
                P \vspace{1mm}\\
                Q \vspace{1mm}\\
                r \vspace{1mm}\\
                s \vspace{1mm}\\
                t \end{pmatrix} = R. $$
Thus $R$ is transformed to $1_-$ by $(\exp\,\varTheta)^{-1} \in (({E_8}^C)^{\sigma,\sigma',\mathfrak{so}(8,C)})_0$.               
\vspace{1mm}

Case (2). $R = (\varPhi, P, Q, r, s, 0), s \not= 0.$ For $\lambda' = \exp(\varTheta\Big(0, 0, 0, 0, \dfrac{\pi}{2}, -\dfrac{\pi}{2}\Big)) \in (({E_8}^C)^{\sigma,\sigma,\mathfrak{so}(8,C)})_0$ (Lemma 4.1.(2)), we have
$$
    \lambda'R = \lambda'(\varPhi, P, Q, r, s, 0) = (\varPhi, Q,-P, -r, 0, -s), \;\;\; -s \not= 0. $$
So, this case can be reduced to Case (1).
\vspace{1mm}

Case (3). $R = (\varPhi, P, Q, r, 0, 0), r \not= 0$. From (2),(5),(6) of Lemma 4.3, we have
$$
      Q \times Q = 0, \;\; \varPhi Q = -3rQ, \;\; \{P, Q\} = 16r^2. $$   
Now, for $\varTheta = \varTheta(0, Q, 0, 0, 0, 0) \in \varTheta(({\mathfrak{e}_8}^C)^{\sigma,\sigma',\mathfrak{so}(8,C)})$ (Lemma 4.1.(2)), we see 
$$
      (\exp\, \varTheta)R = (\varPhi, P + 2rQ, Q, r, -4r^2, 0), \;\; -4r^2 \not= 0. $$
So, this case can be reduced to Case (2).      
\vspace{1mm}

Case (4). $R = (\varPhi, P, Q, 0, 0, 0), Q \not= 0.$ Choose $P_1 \in (\gP^C)_d$ such that $\{P_1, Q\} \not= 0$. For $\varTheta = \varTheta(0, P_1, 0, 0, 0, 0) \in \varTheta(({\mathfrak{e}_8}^C)^{\sigma,\sigma',\mathfrak{so}(8,C)})$ (Lemma 4.1.(2)), we have
$$
       (\exp \, \varTheta)R = \Big(*, *, *, -\dfrac{1}{8}\{P_1, Q\}, *, *\Big). $$
So, this case can be reduced to the Case (3).       
\vspace{1mm}

Case (5). $R = (\varPhi, P, 0, 0, 0, 0), P \not= 0.$ Choose $Q_1 \in (\gP^C)_d$ such that $\{P, Q_1\} \not= 0$. For $\varTheta = \varTheta(0, 0, Q_1, 0, 0, 0) \in \varTheta(({\mathfrak{e}_8}^C)^{\sigma,\sigma',\mathfrak{so}(8,C)})$ (Lemma 4.1.(2)), we have
$$
       (\exp \, \varTheta)R = \Big(*, *, *, \dfrac{1}{8}\{P, Q_1\}, *, *\Big). $$
So, this case can be reduced to the Case (3).       
\vspace{1mm}

Case (6). $R = (\varPhi, 0, 0, 0, 0, 0), \varPhi \not= 0.$ From (10) of Lemma 4.3, we have $\varPhi^2=0$. Choose $P_1 \in (\gP^C)_d$ such that $\varPhi P_1 \not = 0.$  Now, for $\varTheta = \varTheta(0, P_1, 0, 0, 0, 0) \in \varTheta(({\mathfrak{e}_8}^C)^{\sigma,\sigma',\mathfrak{so}(8,C)})$ (Lemma 4.1.(2)), we have
$$
    (\exp\,\varTheta)R = \Big(\varPhi, -\varPhi P_1, 0, 0, \dfrac{1}{8}\{\varPhi P_1, P_1\}, 0 \Big). 
$$
Hence this case is also reduced to Case (2).

Thus the proof of this proposition is completed.  
\end{proof}  
\vspace{3mm}

{\bf Theorem 4.5.} \;\; $({E_8}^C)^{\sigma,\sigma',\mathfrak{so}(8,C)}/(({E_8}^C)^{\sigma,\sigma',\mathfrak{so}(8,C)})_{1_-} \simeq (\gW^C)_{\sigma,\sigma',\mathfrak{so}(8,C)}. $

\noindent {\it In particular, $({E_8}^C)^{\sigma,\sigma',\mathfrak{so}(8,C)}$ is connected.}

\begin{proof} The group $({E_8}^C)^{\sigma,\sigma',\mathfrak{so}(8,C)}$ acts on the space $(\gW^C)_{\sigma,\sigma',\mathfrak{so}(8,C)}$ transitively (Proposition 4.4), 
so the former half of this theorem is proved. The latter half can be shown as follows. Since $(({E_8}^C)^{\sigma,\sigma',\mathfrak{so}(8,C)})_{1_-}$ and $(\gW^C)_{\sigma,\sigma',\mathfrak{so}(8,C)}$ are connected (Propositions 4.2, 4.4),we have that $({E_8}^C)^{\sigma,\sigma',\mathfrak{so}(8,C)}$ is also connected.
\end{proof}
\vspace{3mm}

\section
{\bf  Construction of $Spin(8, C)$ in the group ${E_8}^C$ }
\vspace{3mm}

As similar to the group $({E_8}^C)^{\sigma,\sigma',\mathfrak{so}(8,C)}$, we define the group $(E_8)^{\sigma,\sigma',\mathfrak{so}(8)}$ by
$$
      (E_8)^{\sigma,\sigma',\mathfrak{so}(8)} = \left\{\alpha \in E_8 \, \left| \, \begin{array}{l}
       \sigma\alpha = \alpha\sigma, \sigma'\alpha = \alpha\sigma', \\
       \varTheta(R_{\tilde{D}})\alpha = \alpha\varTheta(R_{\tilde{D}}) \; \mbox{for all} \; \tilde{D} \in \so(8)  \end{array} \right.\right\}, $$
where $\mathfrak{so}(8) = (\mathfrak{f}_4)^{\sigma,\sigma'}$.
\vspace{3mm}
       
{\bf Proposition 5.1.} {\it The Lie algebra $(\mathfrak{e}_8)^{\sigma,\sigma',\mathfrak{so}(8)}$ 
of the group $(E_8)^{\sigma,\sigma',\mathfrak{so}(8)}$ is isomorphic to $\mathfrak{so}(8)$ and the Lie algebra $({\mathfrak{e}_8}^C)^{\sigma,\sigma',\mathfrak{so}(8,C)}$ 
of the group $({E_8}^C)^{\sigma,\sigma',\mathfrak{so}(8,C)}$ is isomorphic to $\mathfrak{so}(8, C)$, that is, we have}
$$
    (\mathfrak{e}_8)^{\sigma,\sigma',\mathfrak{so}(8)} \cong \so(8), \quad ({\mathfrak{e}_8}^C)^{\sigma,\sigma',\mathfrak{so}(8,C)} \vspace{2mm}
\cong \mathfrak{so}(8, C). $$ 

\begin{proof}$(\mathfrak{
e}_8)^{\sigma,\sigma',\mathfrak{so}(8)} \cong \mathfrak{so}(8)$ is proved in [4]. The latter case is the complexification of the former case.
\end{proof}
\vspace{3mm}

{\bf Theorem 5.2.} \qquad \quad $({E_8}^C)^{\sigma,\sigma',\mathfrak{so}(8,C)} \cong Spin(8, C)$.
\begin{proof} The group $({E_8}^C)^{\sigma,\sigma',\mathfrak{so}(8,C)}$ is connected (Theorem 4.5) and its type is $\mathfrak{so}(8, C)$ (Proposition 5.1). Hence the group $({E_8}^C)^{\sigma,\sigma',\mathfrak{so}(8,C)}$ is isomorphic to either one of the following groups
$$
      Spin(8, C), \quad SO(8, C), \quad SO(8, C)/\Z_2. $$      
Their centers of groups above are $\Z_2 \times \Z_2, \Z_2, 1$, respectively. However, we see that the center of $({E_8}^C)^{\sigma,\sigma',\mathfrak{so}(8,C)}$ has $1, \sigma, \sigma', \sigma\sigma'$, so its center is $\Z_2 \times \Z_2$. Hence the group $({E_8}^C)^{\sigma,\sigma',\mathfrak{so}(8,C)}$ have to be $Spin(8, C)$.
\end{proof}
\vspace{5mm}

\section
{\bf The structure of the group $({E_8}^C)^{\sigma,\sigma'}$}
\vspace{3mm}

By using the results above, we shall determine the structure of the group $(\!{E_8}^C\!)^{\sigma,\sigma'} $ $= ({E_8}^C)^\sigma \cap ({E_8}^C)^{\sigma'}$.
\vspace{3mm}

{\bf Lemma 6.1.}\,\,{\it The Lie algebras $({\mathfrak{e}_8}^C)^{\sigma, \sigma'}$ of the groups $({E_8}^C)^{\sigma, \sigma'}$ is given by}
$$
\begin{array}{lll}
 & & ({\mathfrak{e}_8}^C)^{\sigma, \sigma'}
\vspace{2mm}\\
&=&\!\! \{ R \in {\mathfrak{e}_8}^C \,|\, 
           \sigma R = R, \,\, \sigma' R = R \}
\vspace{2mm}\\
 & = &\!\!\left \{ (\varPhi,P, Q, r, s, t) \in {\mathfrak{e}_8}^C \,\left |\,
 \begin{array}{l}
               \varPhi \in ({\mathfrak{e}_7}^C)^{\sigma, \sigma'}, P=(X, Y, \xi, \eta),
     \vspace{1mm}\\
      Q=(Z, W, \zeta, \omega), X, Y \text{\, are diagonal forms},
     \vspace{1mm}\\
Z, W \text{\,are diagonal forms}, \xi, \eta, \zeta, \omega, \in C,
     \vspace{1mm}\\
r, s, t \in C 
\end{array} \right. \right \}.
\end{array}
$$
{\it In particular},
$$
   \dim_{C}(({\mathfrak{e}_8}^C)^{\sigma, \sigma'})=37 + 8 \times 2 + 3=56.
$$
\vspace{1mm}

{\bf Proposition 6.2.}\,\,{\it The group $(({E_8}^C)^{\sigma,\sigma'})_{1_-}$ is a semi-direct product of  \par
\noindent  $\exp(\varTheta(((\gP^C)_d) \oplus C_- ))$ and $({E_7}^C)^{\sigma,\sigma'}$}:
$$
  (({E_8}^C)^{\sigma,\sigma'})_{1_-} = \exp(\varTheta(((\gP^C)_d) \oplus C_- ))\cdot({E_7}^C)^{\sigma,\sigma'}. $$

\noindent {\it In particular, $(({E_8}^C)^{\sigma,\sigma'})_{1_-}$ is connected}.

\begin{proof} We can prove this proposition in a similar way to Proposition 4.2.
\end{proof}
\vspace{2mm}

We denote the connected component of $({E_8}^C)^{\sigma,\sigma'}$ containing the unit element by $(({E_8}^C)^{\sigma,\sigma'})_0$.
\vspace{3mm}

We define a subspace $(\gW^C)_{\sigma,\sigma'}$ of $\gW^C$ by
$$
      (\gW^C)_{\sigma,\sigma'} = \{ R \in \gW^C \, | \, \sigma R = R, \sigma'R = R \}. 
$$
\vspace{1mm}

{\bf Lemma 6.3.} {\it For $R = (\varPhi, P, Q, r, s, t) \in {\mathfrak{e}_8}^C$ satisfying $\sigma R = R, \sigma'R = R, R \not= 0, R$ belongs to $(\gW^C)_{\sigma,\sigma'}$ if and only if $R$ satisfies the following conditions}.
\vspace{1mm}

(1) $2s\varPhi - P \times P = 0$ \quad (2) $2t\varPhi + Q \times Q = 0$
\quad (1) $2r\varPhi + P \times Q = 0$
\vspace{1.5mm}

(4) $\varPhi P -3rP -3sQ = 0$ \quad (5) $\varPhi Q + 3rQ -3tP = 0$
\vspace{1.5mm}

(6) $\{P, Q\} - 16(st + r^2) = 0$
\vspace{1.5mm}

(7) $2(\varPhi P \times Q_1 + 2P \times \varPhi Q_1 - rP \times Q_1 - sQ \times Q_1) - \{P, Q_1\}\varPhi = 0$
\vspace{1.5mm}

(8) $2(\varPhi Q \times P_1 + 2Q \times \varPhi P_1 + rQ \times P_1 -  tP \times  P_1) - \{Q, P_1\}\varPhi = 0$
\vspace{1.5mm}

(9) $8((P \times Q_1)Q - stQ_1 - r^2Q_1 - \varPhi^2Q_1 + 2r\varPhi Q_1) + 5\{P, Q_1\}Q - 2\{Q, Q_1\}P 
$ $= 0$
\vspace{1.5mm}

\hspace*{-1.7mm}(10) $8((Q \times P_1)P + stP_1 + r^2P_1 + \varPhi^2P_1 + 2r\varPhi P_1) + 5\{Q, P_1\}Q - 2\{P, Q_1\}Q
= 0$
\vspace{1.5mm}

\hspace*{-1.7mm}(11) $ 18(\ad\,\varPhi)^2\varPhi_1 + Q \times \varPhi_1P - P \times \varPhi_1Q) + B_7(\varPhi, \varPhi_1)\varPhi = 0$
\vspace{1.5mm}

\hspace*{-1.7mm}(12) $ 18(\varPhi_1\varPhi P - 2\varPhi\varPhi_1P - r\varPhi_1P - s\varPhi_1Q) + B_7(\varPhi, \varPhi_1)P = 0 $
\vspace{1.5mm}

\hspace*{-1.7mm}(13) $ 18(\varPhi_1\varPhi Q - 2\varPhi\varPhi_1Q + r\varPhi_1Q - t\varPhi_1P) + B_7(\varPhi, \varPhi_1)Q = 0 $
\vspace{1.5mm}

\noindent ({\it where $B_7$ is the Killing form of the Lie algebra ${\mathfrak{e}_7}^C$}) {\it for all $\varPhi_1 \in {\mathfrak{e}_7}^C, P_1, Q_1 \in \gP^C$}.
\vspace{2mm}

We denote the connected component of $({E_8}^C)^{\sigma,\sigma'}$ containing the unit element by $(({E_8}^C)^{\sigma,\sigma'})_0$.
\vspace{3mm}

{\bf Proposition 6.4.} {\it The group $(({E_8}^C)^{\sigma,\sigma'})_0$ acts on $(\mathfrak{W}^C)_{\sigma,\sigma'}$ transitively.}

\noindent {\it In particular, $(\gW^C)_{\sigma,\sigma'}$ is connected.}

\begin{proof} We can prove this proposition in a similar way to Proposition 4.4 by using Lemma 6.3.
\end{proof}

\vspace{3mm}

{\bf Proposition 6.5.} \qquad $({E_8}^C)^{\sigma,\sigma'}/(({E_8}^C)^{\sigma,\sigma'})_{1_-} \simeq (\gW^C)_{\sigma,\sigma'}. $

\noindent {\it In particular, $({E_8}^C)^{\sigma,\sigma'}$ is connected.}


\begin{proof} The group $({E_8}^C)^{\sigma,\sigma'}$ acts on $(\gW^C)_{\sigma,\sigma'}$ 
transitively (Proposition 6.4) 
, so the former half of this theorem is proved. Since $(({E_8}^C)^{\sigma,\sigma'})_{1_-}$ and $(\gW^C)_{\sigma,\sigma'}$ are connected (Propositions 6.2, 6.4), we have that $({E_8}^C)^{\sigma,\sigma'}$ is also connected.
\end{proof}
\vspace{3mm}

{\bf Theorem 6.6.} $({E_8}^C)^{\sigma,\sigma'} \cong (Spin(8, C) \times Spin(8, C))/(\Z_2 \times \Z_2), \Z_2 \times \Z_2 = \{(1, 1), (\sigma, \sigma)\} \times \{(1, 1), (\sigma', \sigma')\}$. 

\begin{proof} Let $Spin(8, C) = ({E_8}^C)^{\sigma,\sigma',\mathfrak{so}(8,C)}$ (Theorem 5.2) $\subset ({E_8}^C)^{\sigma,\sigma'}$ and $Spin(8,$ $ C) = ({F_4}^C)^{\sigma,\sigma'}$ (Theorem 2.2) $\subset ({E_8}^C)^{\sigma,\sigma'}$. We define a mapping $\varphi : Spin(8, C) \times Spin(8, C) \to ({E_8}^C)^{\sigma,\sigma'}$ by
$$
     \varphi(\alpha, \beta) = \alpha\beta. 
$$
Since $[R_D, R_8] = 0$ for $R_D \in \spin(8, C) = \so(8, C) = ({\mathfrak{e}_8}^C)^{\sigma,\sigma',\mathfrak{so}(8,C)}$ and $R_8 \in \spin(8, C) = ({\mathfrak{f}_4}^C)^{\sigma,\sigma'}$, we see that $\alpha\beta = \beta\alpha$. Hence $\varphi$ is a homomorphism. $\Ker \, \varphi = \Z_2 \times \Z_2$. In fact, $\dim_C(\spin(8, C) \oplus \spin(8, C)) = 28 + 28 = 56 = \dim_C(({\mathfrak{e}_8})^{\sigma,\sigma'})$ (Lemma 6.1), $\Ker \, \varphi$ is discrete. Hence $\Ker \, \varphi$ is contained  in the center $z(Spin(8, C) \times Spin(8, C)) = z(Spin(8, C)) \times z(Spin(8, C)) = \{1, \sigma, \sigma', \sigma\sigma'\} \times \{1, \sigma, \sigma', \sigma\sigma'\}$. Among them, $\varphi$ maps only $\{(1, 1), (\sigma, \sigma), (\sigma', \sigma'), (\sigma\sigma', \sigma\sigma')\}$ to the identity $1$. Hence we have $\Ker \, \varphi = \{(1, 1), (\sigma, \sigma), (\sigma', \sigma'), (\sigma\sigma', \sigma\sigma')\}= \{(1, 1), (\sigma, $ $\sigma)\}\times \{(1, 1), (\sigma',\sigma')\}= \Z_2 \times \Z_2$.
Since $({E_8}^C)^{\sigma,\sigma'}$ is connected (Proposition 6.5) and again from $\dim_C(({\mathfrak{e}_8}^C)^{\sigma,\sigma'}) = 56 = \dim_C(\spin(8, C) \oplus \spin(8, C))$, we see that $\varphi$ is onto. Thus we have the required isomorphism $({E_8}^C)^{\sigma,\sigma'} \cong (Spin(8, C) \times Spin(8, C))/(\Z_2 \times \Z_2), \Z_2 \times \Z_2 = \{(1, 1), (\sigma,$ $ \sigma)\} \times \{(1, 1), (\sigma', \sigma')\}$.
\end{proof}
\vspace{2mm}

\section
{\bf Main theorem}
\vspace{3mm}

By using results above, we will determine the structure of the group $(E_8)^{\sigma,\sigma'}$ which is the main theorem.
\vspace{3mm}

{\bf Theorem 7.1.} $(E_8)^{\sigma,\sigma'} \cong (Spin(8) \times Spin(8))/(\Z_2 \times \Z_2), \Z_2 \times \Z_2 = \{(1, 1), (\sigma, $ $\sigma)\} \times \{(1, 1), (\sigma', \sigma')\}$. 

\begin{proof} For $\delta \in (E_8)^{\sigma,\sigma'}=(({E_8}^C)^{\tau\wti{\lambda}})^{\sigma,\sigma'} =(({E_8}^C)^{\sigma,\sigma'})^{\tau\wti{\lambda}} \subset ({E_8}^C)^{\sigma,\sigma'}$, there exist $\alpha \in Spin(8, C) = ({E_8}^C)^{\sigma,\sigma',\mathfrak{so}(8,C)}$ and $\beta \in Spin(8, C) = ({F_4}^C)^{\sigma,\sigma'}$ such that $\delta = \varphi(\alpha, \beta) = \alpha\beta$ (Theorem 6.6). From the condition 
$\tau\wti{\lambda}\delta\wti{\lambda}\tau = \delta$, that is, $\tau\wti{\lambda}\varphi(\alpha, \beta)\wti{\lambda}\tau = \varphi(\alpha,\beta)$, we have $\varphi(\tau\wti{\lambda}\alpha\wti{\lambda}\tau, \tau\beta\tau) = \varphi(\alpha, \beta)$. Hence 
$$
\begin{array}{l}
     \mbox{(i)} \;\; \left\{\begin{array}{l}
                      \tau\wti{\lambda}\alpha\wti{\lambda}\tau = \alpha
             \vspace{1mm}\\
                      \tau\beta\tau = \beta,
             \end{array} \right.         
\qquad \qquad
\;      \mbox{(ii)} \;\; \left\{\begin{array}{l}
                      \tau\wti{\lambda}\alpha\wti{\lambda}\tau = \sigma\alpha
\vspace{1mm}\\
                      \tau\beta\tau = \sigma\beta,
             \end{array} \right.       
\vspace{2mm}\\
      \mbox{(iii)} \;\, \left\{\begin{array}{l}
                        \tau\wti{\lambda}\alpha\wti{\lambda}\tau = \sigma'\alpha
               \vspace{1mm}\\
                       \tau\beta\tau = \sigma'\beta,
\end{array} \right.
\qquad \quad 
      \mbox{(iv)} \;\; \left\{\begin{array}{l}
                      \tau\wti{\lambda}\alpha\wti{\lambda}\tau = \sigma\sigma'\alpha
               \vspace{1mm}\\
                      \tau\beta\tau = \sigma\sigma'\beta.
\end{array} \right.        
\end{array} $$

Case (i). The group $\{\alpha \in Spin(8, C) \, |\, \tau\wti{\lambda}\alpha\wti{\lambda}\tau = \alpha\} = (Spin(8, C))^{\tau\wti{\lambda}}$ (which is connected) $= (({E_8}^C)^{\sigma,\sigma',\mathfrak{so}(8,C)})^{\tau\wti{\lambda}}$ (Theorem 5.2) $= ({E_8})^{\sigma,\sigma',\mathfrak{so}(8)}$ and its type is $\so(8)$ (Proposition 5.1), so we see that the group $(E_8)^{\sigma,\sigma',\mathfrak{so}(8)}$ is isomorphic to either one of 
$$
       Spin(8), \quad SO(8), \quad SO(8)/\Z_2. $$
Their centers are $\Z_2 \times \Z_2, \Z_2, 1$, respectively. However the center of $(E_8)^{\sigma,\sigma',\mathfrak{so}(8)}$ has $1, \sigma,\sigma',\sigma\sigma'$, so the group $(E_8)^{\sigma,\sigma',\so(8)}$ has to be $Spin(8)$ and its center is $\{1, \sigma, \sigma', \sigma\sigma'\} = \{1, \sigma\} \times \{1, \sigma'\} = \Z_2 \times \Z_2$. From the condition $\tau\beta\tau = \beta$, we have $\beta \in (Spin(8, C))^\tau = (({F_4}^C)^{\sigma,\sigma'})^\tau = (F_4)^{\sigma,\sigma'} = Spin(8)$ ([3] Theorem 1.5). 
Therefore the group of Case (i) 
\vspace{1mm}
is isomorphic to $(Spin(8) 
\times Spin(8))/(\Z_2 \times \Z_2)$.      

Case (ii). This case is impossible. In fact, set $\beta = (\beta_1, \beta_2, \beta_3), \beta_k \in SO(8, C)$ satisfying $(\beta_1x)(\beta_2y) = \ov{\beta_3(\ov{xy})}, x, y \in \mathfrak{C}^C$. From the condition $\tau\beta\tau = \sigma\beta$, we have $(\tau\beta_1, \tau\beta_2, \tau\beta_3) = (\beta_1, -\beta_2, -\beta_3)$. Hence 
$$
     \tau\beta_1 = \beta_1, \quad \tau\beta_2 = -\beta_2, \quad \tau\beta_3 = - \beta_3. $$
From $\tau\beta_1 = \beta_1$, we have $\beta_1 \in SO(8)$, hence $\beta_2$ is also $\beta_2 \in SO(8)$ by the principle of triality. Then $-\beta_2 = \tau\beta_2 = \beta_2$ which is a contradiction.     
\vspace{1mm}

Case (iii) and (iv) are impossible as similar to Case (ii). 
\vspace{1mm}

Thus we have the required isomorphism $(E_8)^{\sigma,\sigma'} \cong (Spin(8) \times Spin(8))/(\Z_2 \times \Z_2), \Z_2 \times \Z_2 = \{(1, 1), (\sigma, \sigma)\} \times \{(1, 1), (\sigma', \sigma')\}$.
\end{proof}

\vspace{5mm}

\begin{flushright}

\begin{tabular}{l}
 Toshikazu Miyashita \\
Ueda-Higashi High School \\
Ueda City, Nagano, 386-8683,\\
Japan \\
e-mail: anarchybin@gmail.com 
\end{tabular}

\end{flushright}


\bigskip

\vspace{10mm}

\begin{thebibliography}{99}
\bibitem{Andreas} Andreas Kollross,  Exceptional $\mathbb{Z}_2 \times \mathbb{Z}_2$-symmetric spaces. Pacific Journal of Math. 242-1(2009), 113-130.
\vspace{1mm}

\bibitem{Imai} T. Imai and I. Yokota,
\newblock Simply connected compact simple Lie group $E_{8(-248)}$ of type $E_8$, J. Math. Kyoto Univ. 21(1981), 741-762. 
\vspace{1mm}

\bibitem{Miyashita} T. Miyashita and I. Yokota, 
\newblock Fixed points subgroup $G^{\sigma,\sigma'}$ by two involutive automorphisms $\sigma, \sigma'$ of compact exceptional Lie groups $G = F_4, E_6$ and $E_7$, Math. J. Toyama Univ. 24(2001), 135-149. 
\vspace{1mm}

\bibitem{Miyashita} T. Miyashita and I. Yokota, 
\newblock An explicit isomorphism between $(\mathfrak{e}_8)^{\sigma,\sigma',\mathfrak{so}(8)}$ 
and $\mathfrak{so}(8) \oplus \mathfrak{so}(8)$ as Lie algebras, Math. J. Toyama Univ. 24(2001), 
151-158. 
\vspace{1mm}

\bibitem{Yokota} I. Yokota, T. Imai and O. Yasukura, 
\newblock On the homogeneous space $E_8/E_7$, J. Math. Kyoto Univ. 23-3(1983), 467-473.
\vspace{1mm}

\bibitem{Yokotaich} I. Yokota, 
\newblock Realization of involutive automorphisms $\sigma$ and $G^\sigma$ of 
exceptional linear Lie groups $G$, Part I, $G = G_2, F_4$ and $E_6$, Tsukuba J. Math. 14(1990), Part II, $G = E_7$, Tsukuba J. Math. 14(1990), 374-404, Part III, $G = E_8$, Tsukuba J. Math. 15(1991), 301-314.
 
\bibitem{Yokotaichro} I. Yokota, 
\newblock Exceptional simple Lie groups (in Japanese), Gendaisuugakusha(Kyoto). 
\end{thebibliography}
\end{document}